\input ./style/arxiv-general.cfg
\documentclass[aos,MSNbibl,noautosecdot,seceqn,dvips]{arximspdf}
\makeatletter
   \@ifpackageloaded{graphicx}{}{\usepackage{graphicx}}
\makeatother
\usepackage{mathbh}

\doi{10.1214/15-AOS1270C}
\volume{43}
\issue{4}
\pubyear{2015}
\firstpage{1444}
\lastpage{1447}
\docsubty{FLA}
\referstodoi{10.1214/14-AOS1270}

\makeatletter
\let\epsilon\varepsilon
\newcommand{\eqref}[1]{(\ref{#1})}
\makeatother

\begin{document}
\begin{frontmatter}
\vspace*{12pt}
\title{Discussion of ``Frequentist coverage of adaptive nonparametric
Bayesian credible sets''}
\runtitle{Discussion}

\begin{aug}
\author[A]{\fnms{Judith}~\snm{Rousseau}\corref{}\ead
[label=e1]{rousseau@ceremade.dauphine.fr}}
\runauthor{J. Rousseau}
\affiliation{Universit\'e Paris-Dauphine and CREST-ENSAE}
\address[A]{Universit\'e Paris-Dauphine and CREST-ENSAE\\
Place du Mar\'{e}chal de Lattre de Tassigny\\
75016 Paris\\
France\\
\printead{e1}}
\end{aug}

%
\received{\smonth{1} \syear{2015}}


%
\end{frontmatter}

\section{\texorpdfstring{Introduction.}{Introduction}}\label{intro}
I would like first to thank the editors for giving me the opportunity
to discuss the paper \textit{Frequentist coverage of
adaptive Bayesian credible sets} by B. Szabo, A. van der Vaart and
Harry van Zanten. This is a very significant contribution
to the literature on nonparametric credible sets. Before discussing
some specific phenomena enlightened by the paper on frequentist
coverage of adaptive credible sets, I would like to emphasize why I
believe this to be a crucial problem, in particular, in large or
infinite-dimensional models. Over the last ten years or so, there has
been a growing literature on frequentist properties of the posterior
distribution in large or infinite-dimensional models. The results
obtained concern mainly posterior concentration rates, initiated by the
seminal paper of
Ghosal, Ghosh and van~der Vaart
\cite{ggv00}. These results have shown that Bayesian nonparametric
approaches often lead to estimates having very good frequentist
properties such
as adative and minimax (or near minimax with possibly a $\log n$
penalty term) convergence rates under standard loss functions. However,
one of
the interesting aspects of Bayesian methods is that they provide much
more than point estimates via the posterior distribution, in
particular, various sorts of measures of uncertainty can be derived. An
important question is then, how can we understand these measures of
uncertainty? This is all the more important when the models are complex
or large, as it becomes impossible to (1) elicit fully a subjective
prior and (2) to assess perfectly the influence of the prior. Hence,
looking at the frequentist properties of measures of uncertainty is a
way to answer---at least partially---these questions.

In \cite{hoffmanrousseauschmidt14} it is observed that when the
posterior distribution has concentration rate~$\epsilon_n$, under some
loss function $\ell(\cdot, \cdot) $, and if there exists an estimate,
say, the posterior mean, whose convergence rate is also $\epsilon_n$
under the loss $\ell(\cdot, \cdot)$, then if $r_n(\gamma)$ is defined
as the posterior $1-\gamma$ quantile of $ \ell( \theta, \hat\theta)$,
we have
%
\begin{equation}
\label{baycoverage}
1 - \gamma= \int P_\theta \bigl( \ell( \theta, \hat\theta)
\leq r_n(\gamma) \bigr)\,d\pi(\theta), \qquad E_\theta \bigl(
\bigl|r_n(\gamma) \bigr|\bigr) = 0(\epsilon_n),
\end{equation}
where $P_\theta$ and $E_\theta$ denote, respectively, the probability and
expectation under the sampling distribution associated to the parameter
$\theta$.
This means that on average the credible region $\{\theta; \ell(
\theta,
\hat\theta) \leq r_n(\gamma) \}$ has the correct frequentist coverage
and if $\epsilon_n$ is the minimax rate, it has also optimal size,
asymptotically. The question is then, can we describe precisely the set
of parameters $\theta$ such that
%
\begin{equation}
\label{coverage}
P_\theta \bigl( \ell( \theta, \hat\theta) \leq
r_n(\gamma) \bigr) \geq1 -\gamma
\end{equation}
or is at least large enough? In their paper, the authors answer a
similar question by allowing $r_n(\gamma)$ to be inflated by a constant
$L$, possibly large, in the special case of a Gaussian white noise
model with an adaptive empirical Bayes Gaussian prior.

\section{\texorpdfstring{On polished tail parameters.}{On polished tail parameters}}

A key notion in this paper is that of polished tail parameters $\theta
\in\ell_2$,
\[
\sum_{j=N}^\infty\theta_j^2
\leq L_0 \sum_{j=N}^{\rho N}
\theta_j^2 \qquad \forall N > N_0,
\]
where $N_0, L_0$ and $\rho$ are fixed. This definition is not intrinsic
to the function $f = \sum_{j=1}^\infty\theta_j \phi_j$ which is to be
estimated since $f$ can lead to polished tail parameter $\theta$ when
represented in some orthonormal basis $(\phi_j)_j$ and nonpolished tail
parameter $\eta$ in some other orthonormal basis. Although this might
seem disturbing at first, I~think this is inherent to the Bayesian
approach. Let $\pi$ be a prior on some functional set containing a
collection of H\"older balls or Sobolev balls which leads to adaptive
minimax posterior concentration rates over this collection. Then, since
there exists no adaptive concentration rate in $L_2$ and probably in
many other metrics on $f$, the set of \textit{good parameters} (for
which credible regions have good frequentist coverage) is necessarily a
subset of $\Theta$. The question is which subset?
Bull and Nickl
\cite{bullnickl}
(among others) construct procedures such that the subset of \textit{badly behaved} parameters (those that either do not lead to good
coverage or do not lead to optimal size) is rendered as small as
possible. But these procedures require to artificially withdraw the
badly behaved points from the confidence regions, which is not entirely
satisfying. In their case, however, the definition of the well-behaved
parameters does not depend on the representation of the function $f$ in
some particular basis. In the paper of Szabo et al. the advantage is
that the credible region is constructed using standard methods and it
has good frequentist properties for a subset of parameters. Hence, a
natural question arises: Is it possible to construct a prior whose set
of well-behaved parameters is similar to that of Bull and Nickl \cite{bullnickl},
without having to modify the credible region by taking out the badly
behaved parameters?

\section{\texorpdfstring{On the generalization of the results.}{On the generalization of the results}}
 \label{gene}

More importantly, I think that even though the results are presented in
a very specific context, they pave the way for controlling frequentist
coverage of posterior credible regions. Indeed, consider an (inflated)
credible region in the form
\[
C_\alpha^\pi= \bigl\{ \ell( \theta, \hat\theta) \leq
Lr_n(\gamma)\bigr\},
\]
where $r_n(\gamma)$ is the posterior $1-\gamma$ quantile of $\ell(
\theta, \hat\theta)$ and $ \hat\theta$ is some minimax (adaptive)
estimator with rate $\epsilon_n(\theta)$ (the dependence on $\theta$ is
here to emphasize the adaptation property). Under the condition of
posterior concentration rate $\epsilon_n(\theta)$ at~$\theta$,
\[
E_\theta\bigl[r_n(\gamma) \bigr] \lesssim
\epsilon_n;
\]
see \cite{hoffmanrousseauschmidt14} for details on this result. So
the only thing that remains to be verified is that
%
\begin{equation}
\label{minicoverage}
\liminf_n\inf_{\Theta^o}
P_\theta\bigl[\theta\in C_\alpha^\pi\bigr] \geq 1-
\gamma
\end{equation}
for some well-identified subset $\Theta^o $ of $\Theta$.
If the posterior distribution satisfies
%
\begin{equation}
\label{variance}
\Pi \bigl( \theta; \ell( \theta, \hat\theta) > \delta \epsilon
_n(\theta)| X \bigr) \geq\alpha+ o_{P_\theta}(1)
\end{equation}
uniformly on $\Theta^o$, then $P_\theta( r_n(\gamma)> \delta
\epsilon
_n(\theta)) = 1 + o(1)$ uniformly on $\Theta^o $ and
\[
P_\theta \bigl( \theta\in C_\alpha^\pi \bigr) \geq
P_\theta \bigl( \ell( \theta, \hat\theta) \leq L \delta
\epsilon_n(\theta) \bigr) + o(1),
\]
and we can choose $L = 1/\delta$ to ensure \eqref{minicoverage}. So
the main difficulty is to verify \eqref{variance}. When $\hat\theta$ is
the posterior mean and $\ell( \cdot, \cdot) $ is the $L_2$ loss, as in
the present paper, this boils down to bound from below the trace of the
posterior variance. This typically requires that the posterior
distribution asymptotically lives in a space that has an effective
dimension large enough so that the bias of $\hat\theta$ is of the same
order as its variance. The polished tail condition ensures that when
the prior is based on a sequence parameter $\theta\in\ell_2$.

Obviously, this is a rough description of the underlying mechanisms and
this does not take into account some more subtle aspects of the paper.
For instance, the maximum marginal likelihood empirical Bayes approach
is used to simplify the computations since the empirical Bayes
distribution remains Gaussian, while leading to an adaptive posterior
distribution. From the above comments it seems that hierarchical
posterior distribution could be treated using similar ideas, though
less directly. However, how influential are some of the specific
aspects of this empirical Bayes procedure? If instead of maximizing the
marginal likelihood in $\alpha$, the posterior had been computed by
considering a family of priors in the form
\[
\theta_i \stackrel{\mathrm{ind}} {\sim} \mathcal N\bigl(0, \tau
i^{-1-2\alpha} \bigr)
\]
and either consider a hierarchical procedure with a prior on $\tau$ or
a empirical Bayes procedure maximizing the marginal likelihood in $\tau
$, then adaptive posterior concentration rates would be achieved on the
range $\beta\in(0,\alpha+1/2)$ if $\beta$ represents the Sobolev
smoothness of the true parameter. Similar results should be obtained in
this case on the range $\beta\in(0,\alpha+1/2)$. Now, if instead the
prior model had the form $f$, conditional on $\tau$ is a Gaussian
process with kernel
\[
K_\tau(x,y) = e^{-\tau^2 ( x-y)^2}
\]
and $\tau$ follows a Gamma random variable as in \cite{vvvz09}. Then
the posterior has an adaptive concentration rate over collections of H\"
older balls with smoothness $\beta$, with $\beta\in(0, +\infty)$. How
does it impact the behavior of credible regions?

\section{How honest should a confidence region be?}

As I said in Section~\ref{intro}, the questions answered by the authors
in this paper are important questions, as they help to understand some
subtle effects of the prior in large dimensional models. However, as
the nonexistence of adaptive confidence regions over a wide collection
of Sobolev or H\"older classes of functions show, the \textit{full}
minimax paradigm (i.e., having a uniform lower bound on the confidence
and an adaptive minimax upper bound on the size of the confidence
region) has its limits. One might wonder what is the most important?
Weakening the requirement on the confidence or on the size of the
credible regions or considering smaller classes of functions? Somehow
the adaptive Bayesian approach naturally adapts on the size while
losing slightly on the confidence properties of the credible regions,
as shown by \eqref{baycoverage}. Confidence regions constructed in the
frequentist literature are typically honest, however, their sizes are
not uniformly optimal. Using this starting point as a construction of
honest with optimal size confidence regions over smaller functional
classes requires withdrawing from these regions badly behaved
functions. This leads to a somewhat artificial construction. It seems
thus better to be slightly dishonest and start with confidence regions
that have optimal size and to understand over which subclasses of
functions they are honest confidence regions. Obviously, this
construction need not be necessarily Bayesian; however, I believe that
the Bayesian methodology naturally leads to such a construction.






\printaddresses

\begin{thebibliography}{4}


\bibitem{bullnickl}
\begin{barticle}[mr]
\bauthor{\bsnm{Bull},~\bfnm{Adam~D.}\binits{A.~D.}} \AND
\bauthor{\bsnm{Nickl},~\bfnm{Richard}\binits{R.}}
(\byear{2013}).
\btitle{Adaptive confidence sets in {$L\sp 2$}}.
\bjournal{Probab. Theory Related Fields}
\bvolume{156}
\bpages{889--919}.
\bid{doi={10.1007/s00440-012-0446-z}, issn={0178-8051}, mr={3078289}}
\end{barticle}
%

\bptok{imsref}%
\endbibitem

\bibitem{ggv00}
\begin{barticle}[mr]
\bauthor{\bsnm{Ghosal},~\bfnm{Subhashis}\binits{S.}},
\bauthor{\bsnm{Ghosh},~\bfnm{Jayanta~K.}\binits{J.~K.}} \AND
\bauthor{\bsnm{van~der Vaart},~\bfnm{Aad~W.}\binits{A.~W.}}
(\byear{2000}).
\btitle{Convergence rates of posterior distributions}.
\bjournal{Ann. Statist.}
\bvolume{28}
\bpages{500--531}.
\bid{doi={10.1214/aos/1016218228}, issn={0090-5364}, mr={1790007}}
\end{barticle}
%

\bptok{imsref}%
\endbibitem

\bibitem{hoffmanrousseauschmidt14}
\begin{bmisc}[auto:parserefs-M02]
\bauthor{\bsnm{Hoffmann},~\bfnm{M.}\binits{M.}},
\bauthor{\bsnm{Rousseau},~\bfnm{J.}\binits{J.}} \AND
\bauthor{\bsnm{Hieber},~\bfnm{J.~S.}\binits{J.~S.}}
(\byear{2014}).
\bhowpublished{On adaptive posterior concentration rate.
Technical report.}
\end{bmisc}
%

\bptok{imsref}%
\endbibitem

\bibitem{vvvz09}
\begin{barticle}[mr]
\bauthor{\bsnm{van~der Vaart},~\bfnm{A.~W.}\binits{A.~W.}} \AND
\bauthor{\bparticle{van} \bsnm{Zanten},~\bfnm{J.~H.}\binits{J.~H.}}
(\byear{2009}).
\btitle{Adaptive {B}ayesian estimation using a {G}aussian random field with inverse gamma bandwidth}.
\bjournal{Ann. Statist.}
\bvolume{37}
\bpages{2655--2675}.
\bid{doi={10.1214/08-AOS678}, issn={0090-5364}, mr={2541442}}
\end{barticle}
%
\bptok{imsref}%
\endbibitem
\end{thebibliography}
\end{document}